\documentclass{article}
	\usepackage{amsmath,amssymb,amsthm}
	\usepackage{url}

\title{A Proof of Irrationality of $\pi$ by Contraposition\thanks{This research is partially supported by
	JSPS KAKENHI Grant Number 24540071.}}
\author{Akira \sc{Ushijima}\thanks{\texttt{ushijima@se.kanazawa-u.ac.jp}; Faculty of Mathematics and Physics,
Institute of Science and Engineering,
Kanazawa University,
Ishikawa 920-1192,
JAPAN}}
\date{November 25, 2015}

\begin{document}
\maketitle

\begin{abstract}
	We point out that
	the proof of irrationality of $\pi$ by Niven
	can be modified to a proof by contraposition.
	As a warm-up,
	we also give a proof of irrationality of $\sqrt{2}$
	and $\sqrt{3}$
	in a similar way.
\end{abstract}

\section{Introduction}
	Though contradiction and contraposition
	are logically equivalent,
	it is not always clear
	how to translate a proof by contradiction
	to the one by contraposition.
	In this note
	we give a proof of irrationality of $\pi$ by contraposition,
	which is obtained
	by modifying Niven's famous proof \cite{Niven}.
	The strategy of our proof
	is inspired by
	\cite[Theorem~2.3]{Tasaki},
	an irrationality proof of $\sqrt{2}$
	by contradiction,
	with set-theoretic translation.
	Our proof of irrationality of $\pi$
	is given in Section~\ref{sec_pi},
	which follows the section
	that is devoted to the irrationality of $\sqrt{2}$ and $\sqrt{3}$
	to a warm-up of our strategy.

\section{Proof of Irrationality of $\sqrt{2}$ and $\sqrt{3}$} \label{sec_sqrt2}	
	For a given $n \in \mathbb{N}$,
	consider the set $P_n$
	of rational numbers $r$
	with $r^2 = n$.
	Then 
	$n$ is either the square of an odd number or divided by $4$
	if $P_n$ is not the empty set $\varnothing$,
	since the square of any element $r \in P_n$
	is expressed as $4^{p-q} \left( k / l \right)^2$,
	where $r = \left( 2^{p} \, k \right)/ \left( 2^{q} \, l \right)$
	for some odd numbers $k$ and $l$.
	This implies $P_2 = \varnothing$ and $P_3 = \varnothing$.

\section{Proof of Irrationality of $\pi$} \label{sec_pi}
	The main idea of the proof in Section~\ref{sec_sqrt2}
	is to take an indexed family of sets that
	are given as solution sets of equations,
	and then show that
	sets with particular indices are empty by contraposition.
	We apply this idea to Niven's proof \cite{Niven}
	of irrationality of $\pi$.

	Since we use the differentiation of trigonometric functions,
	$\pi$ is supposed to be defined analytically;
	for example,
	it is the smallest positive number for $\sin x = 0$,
	where the trigonometric function $\sin x$ is defined
	as the infinite series
	\begin{equation*}
		\sin x := \sum_{n=0}^{\infty} \frac{\left( -1 \right)^{n} x^{2 \, n + 1}}{\left( 2 \, n + 1 \right) !} .
	\end{equation*}

	For a given rational number
	$r = a/b \in \mathbb{Q}$,
	where $a \in \mathbb{Z}$ and $b \in \mathbb{N}$,
	and a positive integer $n \in \mathbb{N}$,
	define a function $F_n$ by
	\begin{equation*}
		F_n (x) := \sum_{j = 0}^n \left( -1 \right)^j f_n^{( 2 \, j )} (x) ,
	\end{equation*}
	where $f_n (x) := \left( b^n / n! \right) \, x^n \left( r - x \right)^n$.
	The following properties are satisfied for any $r \in \mathbb{Q}$:
	\begin{enumerate}
			\renewcommand{\labelenumi}{(\arabic{enumi})}
		\item \label{prop_equal0r}
			$F_n (0) = F_n (r)$
			for any $n \in \mathbb{N}$;
		\item \label{prop_integer}
			$F_n (0) \in \mathbb{Z}$
			for any $n \in \mathbb{N}$;
		\item \label{prop_nonzero}
			for any $N \in \mathbb{N}$,
			there exists $n \in \mathbb{N}$ such that
			$n \geqq N$ and 
			$F_n (0) \ne 0$.
	\end{enumerate}
	Properties \eqref{prop_equal0r} and \eqref{prop_integer} are
	substantially shown in \cite{Niven}.
	The proof of \eqref{prop_nonzero} is given later.

	Now, let $Q_k$ be the set of positive rational numbers $r$
	satisfying $\sin r = 0$ and $\cos r = k$.
	By the equality $\sin^2 r + \cos^2 r = 1$,
	the number $k$ is either $1$ or $-1$.
	What we need to show is that $\pi \not \in Q_{-1}$.
	It is enough to show that $Q_{-1} = \varnothing$,
	which 
	is done
	by its contrapositive,
	i.e.,
	from now on
	we show that $k$ must be $1$
	if $Q_k \ne \varnothing$.
	
	Since $Q_k \ne \varnothing$,
	we can take an element $r = a / b \in Q_k$.
	For such $r$ and for any $n \in \mathbb{N}$,
	we have
	\begin{align*}
			\int _0^r f_n (x) \, \sin x \, dx
			&= \left[ F_n^{(1)} (x) \, \sin x - F_n (x) \, \cos x \right]_0^r \\
			&= \left( 1 - k \right) F_n (0)
	\end{align*}
	by \eqref{prop_equal0r}.
	
	Since $-1 \leqq \sin x \leqq 1$ and $0 < f_n (x) < \left( b \, r^2 \right)^n / n!$
	for any $x \in (0,r)$ and for any $n \in \mathbb{N}$,
	we have
	\begin{equation*}
		- \frac{\left( b \, r^2 \right)^n}{n!} < f_n (x) \, \sin x < \frac{\left( b \, r^2 \right)^n}{n!}
	\end{equation*}
	on $(0,r)$.
	Since $r > 0$,
	integrating it on $[0,r]$ and we have
	\begin{equation*}
		- \frac{r \left( b \, r^2 \right)^n}{n!} 
			= - \int_0^r \frac{\left( b \, r^2 \right)^n}{n!} dx
			< \int _0^r f_n (x) \, \sin x \, dx
			< \int_0^r \frac{\left( b \, r^2 \right)^n}{n!} dx = \frac{r \left( b \, r^2 \right)^n}{n!} .
	\end{equation*}
	
	This inequality implies
	\begin{equation*}
		- 1 < \int _0^r f_n (x) \, \sin x \, dx = \left( 1 - k \right) F_n (0) < 1
	\end{equation*}
	for any sufficiently large $n \in \mathbb{N}$.
	Since
	$\left( 1 - k \right) F_n (0)$ is an integer
	by \eqref{prop_integer},
	it
	must be $0$.
	Take $n$ to be the one given in \eqref{prop_nonzero}
	and we have $k=1$.
	
	Lastly,
	we prove
	Property
	\eqref{prop_nonzero}.
	For a given $N \in \mathbb{N}$,
	take an odd prime number $p$ with $p > N$,
	and take $n$ to be $p - 1$.
	
	Since $n$ is even,
	the expression of
	$F_n (0)$
	is
	$b^n$ times,
	unlike odd $n$,
	a monic polynomial with integer coefficients up to sign as follows:
	\begin{equation*}
		F_n (0) = b^n \sum_{i=0}^n a_{n,i} \, r^i ,
	\end{equation*}
	where
	\begin{equation*}
		a_{n,i} := \begin{cases}
		0 , & \text{ if $i$ is odd},\\
		\left( -1 \right)^{\frac{i}{2}} \dbinom{n}{n-i} \dfrac{(2 \, n - i)!}{n!} , &
		\text{ if $i$ is even}.
		\end{cases}		
	\end{equation*}
	
	Showing the irreducibility of
	$F_n (0) / b^n$
	in the polynomial ring in $r$ over $\mathbb{Q}$
	is sufficient to prove
	\eqref{prop_nonzero}.
	By Eisenstein's irreducibility criterion,
	it is enough to see the following facts:
	\begin{enumerate}
			\renewcommand{\labelenumi}{(\roman{enumi})}
		\item
			$p$ divides each $a_{n,i}$ for $i = 0 , 1 , \dotsc , n-1$;
		\item
			$p^2$ does not divide $a_{n,0}$.
	\end{enumerate}
	Since
	$a_{n,i} = 0$ for odd $i$,
	we assume that $i = 0 , 2 , \dotsc , n-2$ in the following.
	For such $i$,
	the integer $p = n+1$
	divides $(2 \, n - i)! / n!$.
	We have thus proved (i).
	For (ii),
	since $2 \,n < 2 \, n + 2 = 2 \, p$,
	only the term $n+1$
	in the expression
	$a_{n,0} = (2 \, n)! / n! = 2 \, n \left( 2 \, n - 1 \right) \dotsb \left( n + 1 \right)$
	as a product of integers
	is divided by the prime number $p = n+1$.
	\qed

\end{document}